\documentclass[reqno]{amsart}
\usepackage{amssymb}
\usepackage{amsmath}
\usepackage[mathscr]{euscript}
\usepackage{url}

\newtheorem{lemma}{Lemma}
\newtheorem{theorem}{Theorem}

\newtheorem{proposition}{Proposition}
\newtheorem{definition}{Definition}
\newtheorem{remark}{Remark}

\usepackage{hyperref}
\RequirePackage[dvipsnames]{xcolor} % [dvipsnames]
\definecolor{halfgray}{gray}{0.55} 
\definecolor{webgreen}{rgb}{0,0.5,0}
\definecolor{webbrown}{rgb}{.6,0,0} \hypersetup{%
  colorlinks=true, linktocpage=true, pdfstartpage=3,
  pdfstartview=FitV,%
  breaklinks=true, pdfpagemode=UseNone, pageanchor=true,
  pdfpagemode=UseOutlines,%
  plainpages=false, bookmarksnumbered, bookmarksopen=true,
  bookmarksopenlevel=1,%
  hypertexnames=true,
  pdfhighlight=/O,%hyperfootnotes=true,%nesting=true,%frenchlinks,%
  urlcolor=webbrown, linkcolor=RoyalBlue,
  citecolor=webgreen, %pagecolor=RoyalBlue,%
  pdftitle={},%
  pdfauthor={},%
  pdfsubject={2000 MAthematical Subject Classification: Primary:},%
  pdfkeywords={},%
  pdfcreator={pdfLaTeX},%
  pdfproducer={LaTeX with hyperref}%
}

\begin{document}

\title[Rigidity of fiber bunched cocycles]{Rigidity of fiber bunched cocycles}

\author{Lucas H. Backes}

\address{\noindent IMPA, Estrada Dona Castorina 110, CEP 22460-320 Rio de
  Janeiro, Brasil . 
\newline e-mail: \rm
  \texttt{lhbackes@impa.br} }


\begin{abstract}
We prove a rigidity theorem for fiber bunched matrix-valued H\"{o}lder cocycles over hyperbolic homeomorphisms. More precisely, we show that two such cocycles are cohomologous if and only if they have conjugated periodic data.
\end{abstract}

\maketitle

\section{Introduction}

Given a hyperbolic map $f:M\rightarrow M$ and two $G$-valued H\"{o}lder cocycles $A$ and $B$ where $G$ is a metric group, we are interested in giving necessary and sufficient conditions under which $A$ and $B$ are \textit{cohomologous}, i.e. there exists a $G$-valued H\"{o}lder map $P$ such that  
\begin{center}
$A(x)=P(f(x))B(x)P(x)^{-1}$
\end{center}
for all $x\in M$.\

This problem was initiated by Liv\v{s}ic in his seminal papers \cite{6} and \cite{7} (see also \cite{4}) where he proves that when $G$ has a bi-invariant metric, $A$ is cohomologous to $Id$ i.e. $A$ is a coboundary, if and only if $A(f^{n-1}(p))\ldots A(f(p))A(p)=Id$ for all $p\in M$ such that $f^n(p)=p$ where $Id$ is the identity element of $G$. An important recent development was the extension, by Kalinin \cite{3}, of the Liv\v{s}ic theorem for matrix-valued cocycles.\

With this in mind, it is natural to ask whether two H\"{o}lder cocycles $A$ and $B$ satisfying $A^n(p)=B^n(p)$ for all $p\in M$ such that $f^n(p)=p$ are cohomologous. It is easy to see that this problem reduces to the previous one (i.e. where $B=Id$) when $G$ is abelian, but in general it is not possible to do such a reduction.\

In the present paper we prove that the answer is affirmative, for $G=GL(d,\mathbb{R})$ whenever $A$ and $B$ satisfy a natural fiber bunching condition. The latter is used to construct invariant holonomies, a notion we borrow from Bonatti, G\'{o}mez-Mont, Viana \cite{2,10} and which plays an important part in our arguments.\

A partial case of our result, requiring a strong bounded distortion condition, is contained in the paper of Schmidt \cite{9}. A similar result was also obtained by Parry \cite{8} assuming the group $G$ admits a bi-invariant metric.\

We also point out that the methods we develop in this paper can be extended to cocycles on more general groups, namely, the group of diffeomorphisms of a compact manifold. This is the subject of a joint work with A. Kocsard \cite{BK}.\

After this work had been presented at the 29th Brazilian Mathematical Colloquium \cite{2.1} we learned that similar results were obtained, independently, by Sadovskaya \cite{8.1}.

\section{Definitions and Statements}

Let $(M,d)$ be a compact metric space, $f:M\rightarrow M$ a continuous map and $x\in M$. We define the stable set of $x$ with relation to $f$ as 
\begin{center}$W^s(x)=\lbrace y\in M$; $d(f^n(x),f^n(y))\rightarrow 0$ as $n\rightarrow +\infty \rbrace$
\end{center}
and the stable set of size $\varepsilon >0$ as
\begin{center}$W^s_{\varepsilon}(x)=\lbrace y\in M$; $d(f^n(x),f^n(y))\leq \varepsilon$ for all $n\geq 0 \rbrace$.
\end{center}
Moreover, if $f$ is invertible we can define the unstable set of $x$ and the unstable set of size $\varepsilon >0$ just replacing $f^n$ by $f^{-n}$ in the above definitions.

\begin{definition}
We say that a homeomorphism $f:M\rightarrow M$ is uniformly hyperbolic (or just hyperbolic) if there are constants $C_1$, $\lambda$, $\varepsilon$ and $\tau$ bigger than zero such that
\begin{itemize}
\item $d(f^n(y),f^n(z))\leq C_1e^{-\lambda n}d(y,z)$ for all $y,z\in W^s_{\varepsilon}(x)$ and $n\geq 0$;
\item $d(f^{-n}(y),f^{-n}(z))\leq C_1e^{-\lambda n}d(y,z)$ for all $y,z\in W^u_{\varepsilon}(x)$ and $n\geq 0$;
\item If $d(x,y)<\tau$ then $W^s_{\varepsilon}(x)\cap W^u_{\varepsilon}(y)$ consist of an unique point which we denote by $[x,y]$. Moreover, we require that $[x,y]$ depends continuously on $x$ and $y$.
\end{itemize} 
\end{definition}

Any Anosov diffeomorphism on a compact Riemannian manifold is a hyperbolic homeomorphism. Another interesting class of examples are shifts: let $f:M\rightarrow M$ be the shift map on $M=X^{\mathbb{Z}}$ where $(X,d_{X})$ is a compact metric space and the metric $d$ on $M$ is given by 
$$d((x_n)_n,(y_n)_n)=\sum _{n\in \mathbb{Z}} e^{-\lambda \mid n\mid}\min\lbrace 1,d_{X}(x_n,y_n)\rbrace$$
for any $\lambda >0$.
\begin{definition}
Let $f:M\rightarrow M$ be a homeomorphim on a compact metric space $(M,d)$ and $A:M\rightarrow GL(d,\mathbb{R})$ a continuous map. The pair $(f,A)$ is called a cocycle. Alternatively, we say that the sequence $\lbrace A^n\rbrace _{n\in \mathbb{Z}}$ is a cocycle over $f$ generated by $A$ where
\begin{displaymath}
A^n(x)=
\left\{
	\begin{array}{ll}
		A(f^{n-1}(x))A(f^{n-2}(x))\ldots A(x)  & \mbox{if } n>0 \\
		Id & \mbox{if } n=0 \\
		(A^{-n}(f^{n}(x)))^{-1}=A(f^{n}(x))^{-1}\ldots A(f^{-1}(x))^{-1}& \mbox{if } n<0 \\
	\end{array}
\right.
\end{displaymath}
for all $x\in M$.
\end{definition}

Throughout $\parallel A \parallel$ denotes the operator norm of a matrix $A$, that is, $\parallel A \parallel =\sup \lbrace \parallel Av \parallel ; \parallel v \parallel =1 \rbrace$.

Let us consider $f:M\rightarrow M$ a hyperbolic homeomorphism on a compact metric space $(M,d)$ and $A:M\rightarrow GL(d,\mathbb{R})$ a $\alpha$-H\"{o}lder continuous map which by this we mean that there exists $C_2>0$ such that $\parallel A(x)-A(y)\parallel $ $\leq C_2 d(x,y)^{\alpha}$ for all $x,y\in M$. We say that $A$ satisfies the \textit{fiber bunching condition} if there exists $C_3>0$ and $\theta >0$ with $\theta <\lambda \cdot \alpha$ such that
\begin{equation}\label{eq: fiber-bunching def}
\parallel A^n(x)\parallel \parallel A^n(x)^{-1}\parallel \leq C_3e^{\theta n}
\end{equation}
for all $x\in M$ and $n\geq 0$. We will also use the notion of transitivity: we say that a map $f:M\rightarrow M$ on a metric space $(M,d)$ is \textit{transitive} if there exists $x\in M$ such that $\lbrace f^n(x)\rbrace _{n\in \mathbb{N}}$ is dense in $M$.

\begin{theorem}\label{theorem: main}
Let $f:M\rightarrow M$ be a bi-Lipschitz continuous transitive hyperbolic homeomorphism on a compact metric space $(M,d)$ and assume that $A$, $B:M\rightarrow GL(d,\mathbb{R})$ are $\alpha$-H\"{o}lder continuous maps satisfying the fiber bunching condition. Then the cocycles generated by $A$ and $B$ are cohomologous, i.e. there exist a $\alpha$-H\"{o}lder continuous map $P:M\rightarrow GL(d,\mathbb{R})$ such that 
\begin{center} $A(x)=P(f(x))B(x)P(x)^{-1}$ \end{center}
for all $x\in M$ if and only if there exist a $\alpha$-H\"{o}lder continuous map $Q:M\rightarrow GL(d,\mathbb{R})$ such that $A^n(p)=Q(p)B^n(p)Q(p)^{-1}$ for all $p\in M$ such that $f^n(p)=p$.

\end{theorem} 

As one can easily see, the cohomology relation is an equivalence one over the space of cocycles. In other words, Theorem \ref{theorem: main} says that the periodic data completely characterizes the equivalence classes.

In what follows, for simplicity of the presentation, we will assume $\alpha =1$; the general case is entirely analogous.

\section{Building Invariant Holonomies}

In this section we introduce the main tool used in our proof: invariant holonomies. The existence and main properties of these objects are given by the following proposition which comes from \cite{10}:

\begin{proposition}\label{prop: invariant holonomies}
There exists a constant $C_4>0$ such that, given $x\in M$ and $y, z\in W^s(x)$ the limit 
\begin{center}
$H^{s,A}_{yz}=\lim _{n\rightarrow \infty}A^n(z)^{-1}A^n(y)$
\end{center}
exists and satisfies
\begin{center}
$H^{s,A}_{yz}=H^{s,A}_{xz}H^{s,A}_{yx}$ and $H^{s,A}_{f(y)f(z)}=A(z)H^{s,A}_{yz}A(y)^{-1}$.
\end{center}
Moreover, if $y,z\in W^s_{\varepsilon}(x)$ then
\begin{center}
$\parallel H^{s,A}_{yz}-Id\parallel \leq C_4 d(y,z)$
\end{center}
If $y, z\in W^u_{\varepsilon}(x)$ an analogous result holds for 
\begin{center}
$H^{u,A}_{yz}=\lim _{n\rightarrow \infty}A^{-n}(z)^{-1}A^{-n}(y)$.
\end{center}
An analogous result also holds for $H^{s,B}_{yz}$ and $H^{u,B}_{yz}$.
\end{proposition}

\begin{definition}
$H^{s} $ and $H^{u}$ are called stable and unstable holonomies, respectively.
\end{definition}

To prove Proposition \ref{prop: invariant holonomies}, let us assume that the constant $C_3$ from the definition of fiber bunching in \eqref{eq: fiber-bunching def} is equal to $1$. In particular, $\parallel A(x)\parallel \parallel A(x)^{-1}\parallel \leq  e^{\theta}$ for every $x\in M$. The general case can be deduced from this one.

We will prove only the assertions about $H^{s,A}_{y z}$ since the others are similar and in order to do that we will need the following lemma:

\begin{lemma}\label{lemma: auxiliary 1}
There exists $C=C(A,f)>0$ such that
$$\parallel A^n(y)\parallel \parallel A^n(z)^{-1}\parallel \leq C e^{n\theta}$$
for all $y,z\in W^s_{\varepsilon}(x)$, $x\in M$ and $n\geq 0$.
\end{lemma}

\begin{proof} Note that, since we are assuming $A$ to be Lipschitz, there exists some $L_1=L_1(A)>0$ such that
$$\dfrac{\parallel A(f^j(y))\parallel}{\parallel A(f^j(x))\parallel }\leq \exp (L_1 d(f^j(y),f^j(z)))\leq \exp (L_1 C_1  e^{-j\lambda }2\varepsilon). $$
and similarly for $\dfrac{\parallel A(f^j(z))^{-1}\parallel}{\parallel A(f^j(x))^{-1}\parallel}$. Then, it follows that 
\begin{eqnarray*}
&&\parallel A^n(y)\parallel \parallel A^n(z)^{-1}\parallel \leq \prod ^{n-1}_{j=0} \parallel A(f^j(y))\parallel \cdot \parallel A(f^j(z))^{-1}\parallel \\
&=&\prod ^{n-1}_{j=0} \dfrac{\parallel A(f^j(y))\parallel \cdot \parallel A(f^j(x))^{-1}\parallel}{\parallel A(f^j(x))\parallel } \cdot \dfrac{\parallel A(f^j(z))^{-1}\parallel \cdot \parallel A(f^j(x))\parallel}{\parallel A(f^j(x))^{-1}\parallel} \\
&\leq & L_2 \prod ^{n-1}_{j=0} \parallel A(f^j(x))^{-1}\parallel \cdot \parallel A(f^j(x))\parallel \\
\end{eqnarray*}
where $L_2=\exp (2 L_1 C_1\sum ^{\infty}_{j=0}  e^{-j\lambda }2\varepsilon )$.
Now, using the fiber bunching condition on the last term we get the result.
\end{proof} 

\begin{proof}[Proof of Proposition \ref{prop: invariant holonomies}] By taking forward iterates we can assume that $y,z\in W^s_{\varepsilon}(x)$. For each $n\in \mathbb{N}$ we have that 
$$\parallel A^{n+1}(z)^{-1}A^{n+1}(y)-A^n(z)^{-1}A^n(y)\parallel $$
is bounded by
$$ \parallel A^n(z)^{-1}\parallel \cdot \parallel A(f^n(z))^{-1}A(f^n(y))-Id\parallel \cdot \parallel A^n(y)\parallel. $$
Now, using the Lipschitz continuity of $A$, the contraction property of $f$ on $W^s_{\varepsilon}(x)$ and the previous lemma we get that there exist a constant $L_3=L_3(A,f)>0$ such that
$$\parallel A^{n+1}(z)^{-1}A^{n+1}(y)-A^n(z)^{-1}A^n(y)\parallel \leq L_3e^{n(\theta -\lambda)}d(y,z) .$$
Since $\theta -\lambda <0$, this proves that the sequence is Cauchy and the limit $H^{s,A}_{yz}$ satisfies 
\begin{center}
$\parallel H^{s,A}_{yz}-Id\parallel \leq C_4 d(y,z)$ where $C_4=\sum ^{\infty}_{n=0} L_3e^{n(\theta -\lambda)}$
\end{center} 
as we want.
\end{proof}

Let us fix a constant $C_4>0$ such that Proposition \ref{prop: invariant holonomies} works with this constant for $H^{s,A}_{yz}$, $H^{u,A}_{yz}$, $H^{s,B}_{yz}$, $H^{u,B}_{yz}$ and every (admissible) choice of $y$ and $z$.

\section{Constructing the conjugacy}

Now, with the tool built in the last section in hands, namely, invariant holonomies, we are in position to prove our theorem.

We will prove only one implication since the other one is clear. First of all we note that we can assume without loss of generality that $Q$ is the identity, i.e. that $A^n(p)=B^n(p)$ for all $p\in M$ such that $f^n(p)=p$. This  can be done since we can reduce the previous case to this one just by considering the new cocycle $\tilde{B}(x)=Q(f(x))B(x)Q(x)^{-1}$.

Assume initially that $f$ has a fixed point, i.e. there exists $x\in M$ such that $f(x)=x$ and denote $W(x):=W^s(x)\cap W^u(x)$. It is well known that, in this case, $f$ is topologically mixing, i.e. given any non-empty open sets $U, V\subseteq M$ there exists $k_0 \in \mathbb{N}$ such that $f^{-k}(U)\cap V\neq \emptyset$  for all $k\geq k_0$ [see Corollary 2 of \cite{1}] and $W(x)$ is dense in $M$. Define

$$P:W(x)\rightarrow GL(d,\mathbb{R})$$
by
$$P(y)=H^{s,A}_{xy}(H^{s,B}_{xy})^{-1}=H^{s,A}_{xy}H^{s,B}_{yx}.$$
Note that $P$ satisfies
\begin{center} $A(y)=P(f(y))B(y)P(y)^{-1}$ for all $y\in W(x)$.\end{center}
Indeed,
\begin{eqnarray*}
&&P(f(y))=H^{s,A}_{xf(y)}H^{s,B}_{f(y)x}=H^{s,A}_{f(x)f(y)}H^{s,B}_{f(y)f(x)} \\ &=&A(y)H^{s,A}_{xy}A(x)^{-1}B(x)H^{s,B}_{yx}B(y)^{-1}=A(y)H^{s,A}_{xy}H^{s,B}_{yx}B(y)^{-1}
\end{eqnarray*}
i.e.
\begin{center}
$P(f(y))=A(y)P(y)B(y)^{-1}$
\end{center}
which proves the assertion. To prove it we have used that $f(x)=x$, the hypothesis on the periodic points and Proposition \ref{prop: invariant holonomies}.

Our objective now will be to prove that $P$ is Lipschitz on $W(x)$, so that we can extend it to $M=\overline{W(x)}$. These will be based on a functional identity that we are going to state in Lemma \ref{lemma: transfer map s-u holonomies}. For the proof we will need the following classical fact:

\begin{lemma}\textbf{(Anosov Closing Lemma)} Given $\delta >0$ such that $\theta + \delta <\lambda$ there exist $C_5>0$ and $\varepsilon _0>0$ such that if $z\in M$ satisfy $d(f^n(z),z)<\varepsilon _0$ then there exists a periodic point $p\in M$ such that $f^n(p)=p$ and
 \begin{center}$d(f^j(z),f^j(p))\leq C_5 e^{-(\theta +\delta)\min\lbrace j, n-j\rbrace}d(f^n(z),z) $ for $j=0,1,\ldots ,n$.\end{center}
\end{lemma}

The proof is entirely analogous to the case when the hyperbolic homeomorphism $f$ is, actually, an Anosov diffeomorphism. See for example, \cite{4} p.269, Corollary 6.4.17. 

\begin{lemma}\label{lemma: transfer map s-u holonomies}
$P(y):=H^{s,A}_{xy}H^{s,B}_{yx}=H^{u,A}_{xy}H^{u,B}_{yx}$ for all $y\in W(x)$.
\end{lemma}

\begin{proof}
Fix $y\in W(x)$ and $\delta >0$, $C_5>0$ and $\varepsilon _0>0$ such that Anosov Closing Lemma holds. We begin by noticing that, as $y\in W(x)$ there exist $C_6>0$ and $n_0 \in \mathbb{N}$ such that for all $n\geq n_0$ we have
\begin{center}
$d(f^{-n}(y),f^n(y))\leq C_6e^{-\lambda(n-n_0)}$.
\end{center}
This follows easily from the fact that, as $y\in W(x)=W^s(x)\cap W^u(x)$, there exists $n_0\in \mathbb{N}$ such that $f^{n_0}(y)\in W^s_{\varepsilon }(x)$ and $f^{-n_0}(y)\in W^u_{\varepsilon }(x)$ and that in $W^s_{\varepsilon}(x)$ and $W^u_{\varepsilon}(x)$ we have exponential convergence to $x$. Let $n_1\geq n_0$ be such that for all $n\geq n_1$ we have that $d(f^n(y), f^{-n}(y))<\varepsilon _0$. So, for all $n\geq n_1$ there exists a periodic point $p_n\in M$ with $f^{2n}(p_n)=p_n$ and such that

\begin{center}
$d(f^j(f^{-n}(p_n)), f^j(f^{-n}(y))\leq C_5 e^{-(\theta +\delta)\min \lbrace j, 2n-j\rbrace }d(f^{-n}(y),f^n(y))$ 
\end{center}
for all $j=0,1,\ldots ,2n$. Now, using the hypothesis on the periodic points and that $f^{2n}(f^{-n}(p_n))=f^{-n}(p_n)$ we get that
$$A^{2n}(f^{-n}(p_n))=B^{2n}(f^{-n}(p_n)),$$ 
which can be rewritten as 
$$A^n(p_n)A^n(f^{-n}(p_n))=B^n(p_n)B^n(f^{-n}(p_n)),$$ 
or equivalently as 
$$A^n(f^{-n}(p_n))(B^n(f^{-n}(p_n)))^{-1}=(A^n(p_n))^{-1}B^n(p_n) .$$
Noticing that 
$$A^n(f^{-n}(p_n))= A(f^{-1}(p_n))\ldots A(f^{-n}(p_n))=(A^{-n}(p_n))^{-1}$$
and 
\begin{eqnarray*}
(B^n(f^{-n}(p_n)))^{-1}&=& (B(f^{-1}(p_n))\ldots B(f^{-n}(p_n)))^{-1} \\
&=&B(f^{-n}(p_n))^{-1}\ldots B(f^{-n}(p_n))^{-1}=B^{-n}(p_n)
\end{eqnarray*}
it follows that 
\begin{equation}\label{equa: a}
(A^{-n}(p_n))^{-1}B^{-n}(p_n)=(A^n(p_n))^{-1}B^n(p_n).
\end{equation}

Now we claim that there exist a constant $C_7>0$, independent of $n$ and $p_n$, such that 
$$\parallel (A^n(y))^{-1}B^n(y)-(A^n(p_n))^{-1}B^n(p_n)\parallel<C_7e^{(\theta -\lambda)(n-n_0)}$$
and 
$$\parallel (A^{-n}(y))^{-1}B^{-n}(y)-(A^{-n}(p_n))^{-1}B^{-n}(p_n)\parallel<C_7e^{(\theta -\lambda)(n-n_0)}.$$
Consequently, by \eqref{equa: a} and the fact that $\theta -\lambda <0$ we get that 
$$\parallel (A^n(y))^{-1}B^n(y)-(A^{-n}(y))^{-1}B^{-n}(y)\parallel$$
goes to zero when $n$ goes to infinity. Observing that
$$(A^n(y))^{-1}B^n(y)=(A^n(y))^{-1}A^n(x)B^n(x)^{-1}B^n(y)\longrightarrow _{n\rightarrow \infty} H^{s,A}_{xy}H^{s,B}_{yx}$$
and
$$(A^{-n}(y))^{-1}B^{-n}(y)=(A^{-n}(y))^{-1}A^{-n}(x)B^{-n}(x)^{-1}B^{-n}(y) \longrightarrow _{n\rightarrow \infty} H^{u,A}_{xy}H^{u,B}_{yx}$$
we obtain that
$$P(y)=H^{s,A}_{xy}H^{s,B}_{yx}=H^{u,A}_{xy}H^{u,B}_{yx}$$
as we want. So, to complete the proof it remains to prove our claim. We will only show that there exists $C_7>0$, independent of $n$ and $p_n$, such that
$$\parallel (A^n(y))^{-1}B^n(y)-(A^n(p_n))^{-1}B^n(p_n)\parallel<C_7e^{(\theta -\lambda)(n-n_0)}.$$
The other part is analogous.\

First of all we note that

$$\parallel A^{n}(y)(A^{n}(p_n))^{-1}-Id\parallel $$
is less or equal than
$$ \sum ^{n-1}_{j=0}\parallel A^{n-j}(f^j(y))(A^{n-j}(f^j(p_n)))^{-1}-A^{n-j-1}(f^{j+1}(y))(A^{n-j-1}(f^{j+1}(p_n)))^{-1}\parallel$$
which by the cocycle property is equal to
\begin{eqnarray*}
\sum ^{n-1}_{j=0}\parallel A^{n-j-1}(f^{j+1}(y))A(f^j(y))(A(f^j(p_n)))^{-1}(A^{n-j-1}(f^{j+1}(p_n)))^{-1}\end{eqnarray*}
\begin{flushright}
$- A^{n-j-1}(f^{j+1}(y))(A^{n-j-1}(f^{j+1}(p_n)))^{-1}\parallel .$
\end{flushright}
By the property of the norm this is less or equal than
$$ \sum ^{n-1}_{j=0} \parallel A^{n-j-1}(f^{j+1}(y))\parallel \parallel (A^{n-j-1}(f^{j+1}(p_n)))^{-1}\parallel \parallel A(f^j(y))(A(f^j(p_n)))^{-1} -Id \parallel .$$
Now, using that, since $A$ is Lipschitz, there exist a constant $C_8=C_8(A)>0$ such that 
$$\parallel A(f^j(y))(A(f^j(p_n)))^{-1} -Id \parallel \leq C_8d(f^j(y),f^j(p_n))$$ and that there exist $C_9=C_9(A,\lambda ,\theta ,\delta ,C_5, C_3)>0$ such that, if 
\begin{center}$d(f^j(f^{-n}(y)),f^j(f^{-n}(p_n))\leq C_5 e^{-(\theta +\delta)\min\lbrace j, n-j\rbrace}d(f^{-n}(y),f^n(y)) $ \end{center} 
for $j=0,1,\ldots ,2n$ then
$$\parallel A^{n-j}(f^j(y))\parallel \parallel (A^{n-j}(f^j(p_n))))^{-1}\parallel \leq C_9 e^{\theta (n-j)}$$
for all $j=0,1,\ldots , n$, which proof is similar to the proof of Lemma \ref{lemma: auxiliary 1}, we get that the last expression is less or equal than
$$ \sum ^{n-1}_{j=0} C_8C_9 e^{\theta (n-j-1)}d(f^j(y),f^j(p_n))=\sum ^{n-1}_{j=0} C_8C_9 e^{\theta (n-j-1)}d(f^{n+j}(f^{-n}(y)),f^{n+j}(f^{-n}(p_n)))$$
that by the choice of $p_n$ is less or equal than
$$ \sum ^{n-1}_{j=0} C_5C_8C_9 e^{\theta (n-j-1)}e^{-(\theta +\delta)(n-j)}d(f^{-n}(y),f^n(y))\leq \sum ^{n}_{j=1} C_5C_8C_9 (e^{-\delta})^{j}d(f^{-n}(y),f^n(y))$$
$$\leq \sum ^{\infty}_{j=0} C_5C_8C_9 (e^{-\delta})^{j}d(f^{-n}(y),f^n(y)).$$
Now, defining $C_{10}=C_5C_8C_9\sum ^{\infty}_{j=0} (e^{-\delta})^{j}$ it follows that
$$\parallel A^{n}(y)(A^{n}(p_n))^{-1}-Id\parallel\leq C_{10} d(f^{-n}(y),f^n(y)).$$
Analogously we can prove that
$$\parallel B^{n}(p_n)(B^{n}(y))^{-1}-Id\parallel\leq C_{10} d(f^{-n}(y),f^n(y)).$$
Note now that
\begin{eqnarray*}
&&\parallel A^n(y)(A^n(p_n))^{-1}B^n(p_n)(B^n(y))^{-1} - Id \parallel \\
&\leq &\parallel A^n(y)(A^n(p_n))^{-1}\parallel \parallel B^n(p_n)(B^n(y))^{-1} - Id \parallel + \parallel A^n(y)(A^n(p_n))^{-1} - Id \parallel \\
&\leq & C_{11}d(f^n(y),f^{-n}(y))
\end{eqnarray*}
since $\parallel A^n(y)(A^n(p_n))^{-1}\parallel \leq C_{10}$diam$(M)+1$, independent of $n$ and $p_n$ where $C_{11}=2(C_{10}$diam$(M)+1)$. Consequently, using that $A^n(x)=B^n(x)$, since $f(x)=x$, and that $\parallel A^n(x)\parallel \parallel (A^n(x))^{-1}\parallel \leq C_3e^{\theta n}$ we get that
\begin{eqnarray*}
&&\parallel (A^n(x))^{-1}A^n(y)(A^n(p_n))^{-1}B^n(p_n)(B^n(y))^{-1}B^n(x) - Id \parallel \\
&\leq & \parallel (A^n(x))^{-1}\parallel \parallel A^n(x)\parallel \parallel A^n(y)(A^n(p_n))^{-1}B^n(p_n)(B^n(y))^{-1} - Id \parallel \\
&\leq & C_3 e^{\theta n}\parallel A^n(y)(A^n(p_n))^{-1}B^n(p_n)(B^n(y))^{-1} - Id \parallel \\
&\leq & C_3C_{11}e^{\theta n}d(f^{-n}(y),f^n(y))\leq C_3C_{11}e^{\theta n}C_6e^{-\lambda (n-n_0)}=C_3C_6C_{11}e^{\theta n_0}e^{(\theta -\lambda) (n-n_0)}
\end{eqnarray*}
that is,
$$ \parallel (A^n(x))^{-1}A^n(y)(A^n(p_n))^{-1}B^n(p_n)(B^n(y))^{-1}B^n(x) - Id \parallel \leq C_{12}e^{(\theta -\lambda) (n-n_0)}$$
where $C_{12}=C_3C_6C_{11}e^{\theta n_0}$. Now, as $(A^n(y))^{-1}A^n(x)$ converges to $ H^{s,A}_{xy}$ when $n$ goes to infinity and $(B^n(x))^{-1}B^n(y)$ converges to $H^{s,B}_{yx}$, there exist $N>0$ such that

\begin{center}
$\parallel(A^n(y))^{-1}A^n(x)\parallel \leq N$ and $\parallel(B^n(x))^{-1}B^n(y)\parallel \leq N$
\end{center}
for all $n\in \mathbb{N}$. Thus,
\begin{eqnarray*}
&&\parallel (A^n(p_n))^{-1}B^n(p_n)- (A^n(y))^{-1}B^n(y)\parallel \\
&=& \parallel (A^n(p_n))^{-1}B^n(p_n)- (A^n(y))^{-1}A^n(x)(B^n(x))^{-1}B^n(y)\parallel \\
&\leq &\parallel (A^n(y))^{-1}A^n(x)\parallel \parallel (B^n(x))^{-1}B^n(y)\parallel \\
&\cdot &\parallel (A^n(x))^{-1}A^n(y)(A^n(p_n))^{-1}B^n(p_n)(B^n(y))^{-1}B^n(x)- Id\parallel \\
&\leq &N^2 C_{12}e^{(\theta -\lambda) (n-n_0)} 
\end{eqnarray*}
Therefore, considering $C_7=N^2 C_{12}$ we get that 
$$\parallel (A^n(p_n))^{-1}B^n(p_n)- (A^n(y))^{-1}B^n(y)\parallel \leq C_{7}e^{(\theta -\lambda) (n-n_0)}$$
which proves our claim and consequently the lemma.
\end{proof}

\begin{lemma}\label{lemma: P is Lipschitz}
$P$ is Lipschitz on $W(x)$.
\end{lemma}

\begin{proof}
Let $\varepsilon $ and $ \tau >0$ be given as in the definition of hyperbolic homeomorphism. Take $y,z\in W(x)$ and assume initially that $z\in W^s_{\varepsilon }(y)$. Then, 
\begin{eqnarray*}
&&\parallel P(y)-P(z) \parallel = \parallel H^{s,A}_{xy}H^{s,B}_{yx}-H^{s,A}_{xz}H^{s,B}_{zx} \parallel \\
&=& \parallel (H^{s,A}_{xy}(H^{s,A}_{xz})^{-1}-Id)H^{s,A}_{xz}H^{s,B}_{yx} +  H^{s,A}_{xz}H^{s,B}_{yx}(Id - (H^{s,B}_{yx})^{-1}H^{s,B}_{zx}) \parallel \\
&=& \parallel (H^{s,A}_{zy}-Id)H^{s,A}_{xz}H^{s,B}_{yx} +  H^{s,A}_{xz}H^{s,B}_{yx}(Id - H^{s,B}_{zy}) \parallel \\
&\leq &\parallel H^{s,A}_{zy}-Id\parallel \parallel H^{s,A}_{xz}H^{s,B}_{yx}\parallel + \parallel H^{s,A}_{xz}H^{s,B}_{yx}\parallel \parallel Id - H^{s,B}_{zy} \parallel \\
&\leq &2 C_4\parallel H^{s,A}_{xz}H^{s,B}_{yx}\parallel d(y,z)
\end{eqnarray*}
where the last inequality is due to Proposition \ref{prop: invariant holonomies}. Observe now that
$$\parallel H^{s,A}_{xz}H^{s,B}_{yx}\parallel = \parallel H^{s,A}_{xz} H^{s,A}_{yx}H^{s,A}_{xy} H^{s,B}_{yx}\parallel \leq \parallel H^{s,A}_{yz}\parallel \parallel P(y)\parallel .$$
Thus, as $\parallel H^{s,A}_{yz}\parallel \leq C_4d(y,z)+1 \leq C_4 \varepsilon  +1$ by Proposition \ref{prop: invariant holonomies}, it follows that
$$\parallel P(y)-P(z) \parallel \leq 2C_4(C_4 \varepsilon +1)\parallel P(y)\parallel d(y,z)$$ 
whenever $y,z\in W(x)$ and $z\in W^s_{\varepsilon }(y)$. Using Lemma \ref{lemma: transfer map s-u holonomies} and proceeding in the same way we get that, if $y,z\in W(x)$ and $z\in W^u_{\varepsilon }(y)$ then
$$\parallel P(y)-P(z) \parallel \leq 2C_4(C_4 \varepsilon  +1)\parallel P(y)\parallel d(y,z).$$
Now, given $y,z\in W(x)$ with $d(y,z)<\tau $ let us consider $w=[y,z]=W^s_{\varepsilon }(y)\cap W^u_{\varepsilon }(z)$. Therefore, as $w\in W^s_{\varepsilon }(y)$ and $w\in W^u_{\varepsilon }(z)$ it follows by the previous comments that 
$$\parallel P(y)-P(z) \parallel \leq \parallel P(y)-P(w) \parallel +\parallel P(w)-P(z) \parallel $$
$$\leq  2C_4(C_4 \varepsilon  +1)\parallel P(w)\parallel (d(y,w) + d(w,z)).$$
Then, using that $\parallel P(w)\parallel \leq C_{13} \parallel P(y)\parallel$ where $C_{13}=2C_4(C_4 \varepsilon  +1)+1$ does not depend on $y$ and there exists $C_{14}>0$ such that $d(y,w)+ d(w,z)\leq C_{14}d(y,z)$ where this constant does not depend on $y,w$ or $z$ since $w=[y,z]=W^s_{\varepsilon }(y)\cap W^u_{\varepsilon }(z)$ it follows that
$$\parallel P(y)-P(z) \parallel \leq C_{15}\parallel P(y)\parallel d(y,z)$$
where $C_{15}=2C_4(C_4 \varepsilon  +1)C_{13}C_{14}>0$ does not depend neither on $y$ nor $z$.

Let us consider now $\lbrace y_1,y_2,\ldots ,y_k\rbrace \subset W(x)$ such that $M\subseteq \cup ^k_{j=1}B(y_j,\tau )$ where $B(y_j,\tau )$ is the ball of radius $\tau $ and centred on $y_j$ and 
$$C_{16}\geq \max \lbrace \max _{j=1,\ldots ,k}C_{15}\parallel P(y_j)\parallel \varepsilon  , \max _{j=1,\ldots ,k}\parallel P(y_j)\parallel \rbrace .$$
Note that there exists such $y_j$'s since $W(x)$ is dense in $M$ and $M$ is compact. So, given $y\in W(x)$ there exist $j\in \lbrace 1,2,\ldots ,k \rbrace$ such that $y\in B(y_j,\tau )$. By the previous considerations we get that 
$$\parallel P(y)-P(y_j)\parallel \leq C_{15}\parallel P(y_j)\parallel \varepsilon  \leq C_{16}$$
and consequently $\parallel P(y)\parallel \leq 2C_{16}.$ Thus, given $y,z\in W(x)$ with $d(y,z)<\tau $ it follows that 
$$\parallel P(y)-P(z)\parallel \leq 2C_{15}C_{16}d(y,z)$$ 
and consequently $P$ is Lipschitz as we want.
\end{proof}

\section{Concluding the proof of the main result}

Let us assume initially that $f$ is topologically mixing. If $f$ has a fixed point, by Lemma \ref{lemma: P is Lipschitz} we know that $P$ is Lipschitz and thus we can extend it to $\overline{W(x)} =M$. Such extension $\overline{P}$ is also Lipschitz and satisfies $A(y)=\overline{P}(f(y))B(y)\overline{P}(y)^{-1}$ for all $y\in \overline{W(x)}=M$ since $P$ satisfies it on $W(x)$. Thus, $\overline{P}$ is a Lipschitz map satisfying the desired identity.  For the general case, let $x\in M$ be a periodic point and $n_0$ its period. Consider now the new cocycles $\tilde{A},\tilde{B}:M\rightarrow GL(d,\mathbb{R})$ over $F=f^{n_0}$ given by
$$\tilde{A}(x)=A(f^{n_0 -1}(x))\ldots A(x)$$
and
$$\tilde{B}(x)=B(f^{n_0 -1}(x))\ldots B(x).$$
It is easy to see that $\tilde{A}$ and $\tilde{B}$ satisfy the fiber bunching condition (over $F$), are Lipschitz and $F(x)=x$. Thus, applying the previous result to this case we get that there exists a Lipschitz continuous map $P:M\rightarrow GL(d,\mathbb{R})$ such that 
$$\tilde{A}(y)=P(F(y))\tilde{B}(y)P(y)^{-1}$$
for all $y\in M$. Rewriting this in terms of the original cocycles we get that there exists a Lipschitz continuous map $P:M\rightarrow GL(d,\mathbb{R})$ such that
$$A^{n_0}(y)=P(f^{n_0}(y))B^{n_0}(y)P(y)^{-1}$$
for all $y\in M$.\

Now, in order to complete the proof, let us prove a couple of lemmas:

\begin{lemma}
There exists $P:M\rightarrow GL(\mathbb{R},d)$, a Lipschitz continuous map, such that 
$$A(y)=P(f(y))B(y)P(y)^{-1}$$
for all $y\in M$ if and only if there exist $P:M\rightarrow GL(\mathbb{R}, d)$, a Lipschitz continuous map, and $m,n \in \mathbb{N}$ relatively prime such that 
$$A^n(y)=P(f^n(y))B^n(y)P(y)^{-1}$$
and 
$$A^m(y)=P(f^m(y))B^m(y)P(y)^{-1}$$ for all $y\in M$.

\end{lemma}

\begin{proof}
One implication is trivial. Let us prove the other one. Assume that there exist $P:M\rightarrow GL(\mathbb{R}, d)$ Lipschitz and $m,n \in \mathbb{N}$ relatively prime such that 
$$A^n(y)=P(f^n(y))B^n(y)P(y)^{-1}$$
and 
$$A^m(y)=P(f^m(y))B^m(y)P(y)^{-1}$$ for all $y\in M$. As $m$ and $n$ are relatively prime there exist $k,l \in \mathbb{Z}$ such that $1=mk+nl$. Then,
\begin{eqnarray*}
&&A(y)=A^{mk+nl}(y)=A^{mk}(f^{nl}(y))A^{nl}(y)\\ 
&=&P(f^{mk}(f^{nl}(y)))B^{mk}(f^{nl}(y))P(f^{nl}(y))^{-1}P(f^{nl}(y))B^{nl}(y)P(y)^{-1} \\
&=& P(f(y))B^{mk}(f^{nl}(y))B^{nl}(y)P(y)^{-1}= P(f(y))B(y)P(y)^{-1}
\end{eqnarray*}
for all $y\in M$ proving the lemma.
\end{proof}

\begin{lemma}
There exists $P:M\rightarrow GL(\mathbb{R}, d)$, a Lipschitz continuous map, and $m,n \in \mathbb{N}$ relatively prime such that 
$$A^n(y)=P(f^n(y))B^n(y)P(y)^{-1}$$
and 
$$A^m(y)=P(f^m(y))B^m(y)P(y)^{-1}$$ for all $y\in M$.
\end{lemma}

\begin{proof}
Let $\lbrace p_1, p_2, \ldots ,p_{d^2}, p_{d^2+1}\rbrace \subset M$ be a subset of distinct periodic points such that the periods $n_i$ of $p_i$ and $n_j$ of $p_j$ are relatively prime for all $i\neq j$ and denote $N:=\prod^{d^2+1}_{j=1}n_j$. The existence of such periodic points is ensured by the hypothesis on $f$, namely, hyperbolicity and the topological mixing property. Given $j\in \lbrace 1, 2, \ldots ,d^2+1\rbrace$ consider $P_{p_j}:M\rightarrow GL(\mathbb{R},d)$ the solution of 
$$A^{n_j}(z)=P_{p_j}(f^{n_j}(z))B^{n_j}(z)P_{p_j}(z)^{-1}$$
for all $z\in M$ given by the procedure described above and fix $y\in M$ such that its orbit under $f^N$ is dense in $M$ which exists since $f$ is topologically mixing. Then, as the dimension of $GL(\mathbb{R},d)$ is $d^2$, there exist $p\in \lbrace p_1, p_2, \ldots ,p_{d^2}, p_{d^2+1}\rbrace $, which we will assume without loss of generality to be $p_{d^2+1}$, and $c_j \in \mathbb{R}$, $j=1,2,\ldots ,d^2$ such that 
$$P_p(y)=P_{p_{d^2+1}}(y)=\sum ^{d^2}_{j=1}c_jP_{p_j}(y).$$
Note that this equality holds, a priori, just for this fixed $y$. Now, let $l$ be a multiple of $n_1, n_2,\ldots ,n_{d^2}$ and $n:=n_{d^2+1}$ (period of $p$). Then, 
\begin{eqnarray*}
P_p(f^l(y))&=&A^l(y)P_p(y)B^l(y)^{-1}=A^l(y) \sum ^{d^2}_{j=1}c_jP_{p_j}(y) B^l(y)^{-1} \\
&=&\sum ^{d^2}_{j=1}c_j A^l(y)P_{p_j}(y)B^l(y)^{-1}= \sum ^{d^2}_{j=1}c_j P_{p_j}(f^l(y))\\
\end{eqnarray*}
by the choice of $l$ and the properties of $P_{p_j}$. Thus, as the orbit of $y$ under $f^N$ is dense in $M$ and $P_{p_j}$ is Lipschitz continuous for all $j$, it follows that
$$P_p(z)=\sum ^{d^2}_{j=1}c_jP_{p_j}(z)$$
for all $z\in M$.
Let us consider now $m\in \mathbb{N}$ such that $m$ is multiple of $n_1,n_2, \ldots ,n_{d^2}$ and relatively prime with $n$. Then, 
$$P_p(f^m(z))=A^m(z)P_p(z)B^m(z)^{-1}$$
for all $z\in M$.
Indeed, 
\begin{eqnarray*}
P_p(f^m(z))&=&\sum ^{d^2}_{j=1}c_jP_{p_j}(f^m(z))=\sum ^{d^2}_{j=1}c_jA^m(z)P_{p_j}(z)B^m(z)^{-1} \\
&=& A^m(z)\sum ^{d^2}_{j=1}c_jP_{p_j}(z) B^m(z)^{-1}= A^m(z)P_p(z)B^m(z)^{-1} \\
\end{eqnarray*}
which completes the proof of this lemma since $m$ and $n$ are relative primes and 
$$A^n(z)=P_p(f^n(z))B^n(z)P_p(z)^{-1}$$
and 
$$A^m(z)=P_p(f^m(z))B^m(z)P_p(z)^{-1}$$ for all $z\in M$ as we want.
\end{proof} 

Thus, combining this two lemmas we complete the proof of Theorem \ref{theorem: main} in the case when $f$ is topologically mixing. Assume now that $f$ is transitive. By the Spectral Decomposition Theorem [see Theorem 3 of \cite{1}] there exist pairwise disjoint compact sets $\Omega _1, \Omega _2 , \ldots ,\Omega _k$ such that $M=\Omega _1 \cup \Omega _2 \cup \ldots \cup\Omega _k$, $f(\Omega _i)=\Omega _{i+1}$ for $i=1,\ldots ,k-1$ and $f(\Omega _k)=\Omega _1$ and such that $f^k\mid _{\Omega _i}:\Omega _i \rightarrow \Omega _i$ is topologically mixing for $i=1, \ldots ,k$. Let us consider $f^k\mid _{\Omega _1}:\Omega _1 \rightarrow \Omega _1$ and $A^k,B^k:\Omega _1\rightarrow GL(d,\mathbb{R})$ as usual given by $A^k(x)=A(f^{k-1}(x))\ldots A(x)$ and $B^k(x)=B(f^{k-1}(x))\ldots B(x)$. Thus, applying the previous results we get that there exists $P_1:\Omega _1 \rightarrow GL(d,\mathbb{R})$, a Lipschitz continuous map, such that 
\begin{center}$A^k(x)=P_1(f^k(x))B^k(x)P_1(x)^{-1}$ for all $x\in \Omega _1$ .\end{center}
Define $P:M\rightarrow GL(d,\mathbb{R})$ by
\begin{center}
$P(f^j(x))=A^j(x)P_1(x)B^j(x)^{-1}$ for $x\in \Omega _1$ and $j=0,1,\ldots ,k-1$.
\end{center}
Such $P$ is well defined and Lipschitz continuous and it satisfies 
$$A(y)=P(f(y))B(y)P(y)^{-1}$$
for all $y\in M$. Indeed, given $x\in \Omega _1$ and $j\in \lbrace 0,1,\ldots ,k-2\rbrace$ 
$$P(f^{j+1}(x))=A^{j+1}(x)P_1(x)B^{j+1}(x)^{-1}$$
$$ = A(f^j(x))A^j(x)P_1(x)B^j(x)^{-1}B(f^j(x))^{-1}=A(f^j(x))P(f^{j}(x))B(f^{j}(x))^{-1}$$ 
as we want. If $j=k-1$ we have that 
$$P(f^k(x))=P_1(f^k(x))=A^k(x)P_1(x)B^k(x)^{-1}$$
and
$$P(f^{k-1}(x))=A^{k-1}(x)P_1(x)B^{k-1}(x)^{-1}$$
$$=A(f^{k-1}(x))^{-1}A^{k}(x)P_1(x)B^k(x)^{-1}B(f^{k-1}(x))$$
$$=A(f^{k-1}(x))^{-1}P(f^k(x))B(f^{k-1}(x))$$
which completes the proof Theorem \ref{theorem: main}. $\hfill \square $

\begin{remark}
We could have assumed in the statement of Theorem \ref{theorem: main} that only one of the cocycles satisfies the fiber bunching condition since the other one would satisfy it automatically. This follows from the hypothesis on periodic points and Proposition 2.1 of \cite{4.0}, which says that given an $\alpha$-H\"{o}lder cocycle $F:M\rightarrow GL(d,\mathbb{R})$ over a hyperbolic homeomorphism $f:M\rightarrow M$ such that $\lambda ^+(F,p)-\lambda ^-(F,p)<\theta $ for every $f$-periodic point $p\in M$ and some $\theta <\lambda \alpha$ where $\lambda ^+(F,p)$ denotes the largest Lyapunov exponent of $F$ at $p$ and $\lambda ^-(F,p)$ denotes the smallest one, then $F$ satisfies the fiber bunching condition.

\end{remark}

\medskip{\bf Acknowledgements.} The author would like to express his sincere gratitude to his advisor Professor Marcelo Viana for his guidance throughout the preparation of this work. The author was supported by a CNPq-Brazil doctoral fellowship.

\end{document}